%mourad.tex:
%Efficient Weighted Counting of Multiset Derangements
%%a Plain TeX file by Shalosh B. Ekhad and Doron Zeilberger (x pages)

%begin macros

\baselineskip=14pt
\parskip=10pt

\font\eightrm=cmr8 

\magnification=\magstephalf

\def\1{{\overline{1}}}
\def\2{{\overline{2}}}
\parindent=0pt
\overfullrule=0in

\def\frac#1#2{{#1 \over #2}}
%\headline={\rm  \ifodd\pageno  \RightHead  \else  \LeftHead  \fi}
%\def\RightHead{\centerline{
%Title
%}}
%\def\LeftHead{ \centerline{Doron Zeilberger}}
%end macros
\centerline
{\bf  
Efficient Weighted Counting of Multiset Derangements
}
\bigskip
\centerline
{\it Shalosh B. EKHAD and Doron ZEILBERGER}

\bigskip

\qquad \qquad {\it Dedicated to Mourad El-Houssieny ISMAIL (born April 27, 1944), with friendship and admiration }

{\bf Abstract}: We use the Almkvist-Zeilberger algorithm, combined with a weighted version of the Even-Gillis Laguerre integral due to  Foata and Zeilberger, 
in order to efficiently compute weight enumerators of multiset derangements according to the number of cycles.
The present paper is inspired by important previous work by Mourad Ismail and his collaborators, done in the late 1970s, but still useful
after all these years.

{\bf Multiset Derangements}

The (generalized) Laguerre polynomials, $L_k^{(\alpha)}(x)$, are defined as follows:
$$
L_k^{(\alpha)}(x) \, :=  \, \sum_{i=0}^{n} (-1)^i \, \frac{(\alpha+i+1) (\alpha+i+2)  \dots  (\alpha+n)}{i!(n-i)!}\, x^i \quad .
$$

Let $A_1,A_2, \dots, A_n$ be $n$ pairwise-disjoint sets of cardinalities $k_1, \dots, k_n$ respectively. A permutation $\pi$ of $A:=A_1 \cup \dots \cup A_n$ is a multiset derangement
if for every $1 \leq i \leq n$, whenever $x \in A_i$, $\pi(x) \not \in A_i$. Let $D(k_1, \dots, k_n)$ be the set of such multiset derangements.

As usual, for any finite set $S$, let $|S|$ denote its number of elements.

Shimon Even and Joe Gillis [EvG] (see also [As], [G], and [Z]) proved that 
$$
|D(k_1, \dots, k_n)|=(-1)^{k_1 + \dots + k_n} \left (\prod_{i=1}^{n} k_i! \right) \int_{0}^{\infty} \left ( \prod_{i=1}^{n} L^{(0)}_{k_i}(x) \right ) e^{-x} \, dx \quad.
\eqno(1)
$$

{\bf Comment:} Usually the elements of each $A_i$ are identified, and the formula then does not have  $\left (\prod_{i=1}^{n} k_i! \right)$ in front, but for the purpose of the
present paper, as was done in [FZ], we have $k_1+ \dots +k_n$ distinct elements, all preserving their individuality.

For a permutation $\pi$, let $cyc(\pi)$ denote its {\it number of cycles}. For example $cyc(1234)=4$, and $cyc(2341)=1$. Recall that, famously,
$$
\sum_{\pi \in S_n} \alpha^{cyc(\pi)} \, = \, \alpha(\alpha+1) \cdots (\alpha+n-1) \quad.
$$

{\eightrm
This is easily proved by considering how to form permutations of $\{1, ...,n\}$ out of those of  $\{1, ...,n-1\}$.  Denote the left side by $f_n(\alpha)$.
Given a permutation of  $\{1, ...,n-1\}$ we can insert $n$ inside any of the existing cycles, and there are $n-1$ ways of doing it, 
preserving the number of cycles, or else create a brand-new cycle with $n$ alone, increasing the number of cycles by one. Hence $f_n(\alpha)=((n-1)+ \alpha)f_{n-1}(\alpha)$.
}

Let  $A(k_1, \dots, k_n)(\alpha)$ be the weight enumerator, according to the number of cycles, of the set of multiset derangements, $D(k_1, \dots, k_n)$, in other words
$$
A(k_1, \dots, k_n)(\alpha) \, := \, \sum_{\pi \in D(k_1, \dots, k_n) } \alpha^{cyc(\pi)} \quad .
$$

In 1988, Dominique Foata and one of us (DZ) [FZ], inspired by the work of Mourad Ismail and his collaborators ([AsIs], [AsIsRa], [AsIsKo]) proved the following $\alpha$-analog of
the Even-Gillis theorem:

$$
A(k_1, \dots, k_n)(\alpha)=\frac{(-1)^{k_1 + \dots + k_n}}{(\alpha -1)!} \left (\prod_{i=1}^{n} k_i! \right) \int_{0}^{\infty} \left ( \prod_{i=1}^{n} L^{(\alpha -1)}_{k_i}(x) \right ) x^{\alpha -1}  e^{-x} \, dx \quad.
\eqno(2)
$$

Note that this is a polynomial of degree $(n_1 + \dots +n_k)/2$ rather than  $n_1+ \dots + n_k$, since every cycle is at least of length $2$.

In this paper we will focus on efficient computations  of many terms of the sequences $A(k,k, \dots, k)(\alpha)$, where $k$ is repeated $n$ times,
for specific (small, and not so small) $k$, but arbitrarily large $n$. In other words, our goal  is to compute as many as possible terms of the sequences
$$
F_k(n)(\alpha):=A(k, \dots, k) (\alpha) \quad, \quad (k \quad repeated \quad n \quad times) \quad,
$$
for $k=1$, $k=2$, etc.

By $(2)$ we have
$$
F_{k}(n)(\alpha)= \frac{(-1)^{kn} (kn)!}{(\alpha -1)!} \int_{0}^{\infty} (L^{(\alpha -1)}_{k}(x))^{n}   x^{\alpha -1} e^{-x} \, dx \quad.
\eqno(3)
$$

Using the  {\bf Almkvist-Zeilberger algorithm} [AlZ] (see [D] for a lucid and engaging account), one of us (SBE), using our Maple package {\tt Mourad.txt}, available from the front of this article

{\tt https://sites.math.rutgers.edu/\~{}zeilberg/mamarim/mamarimhtml/mourad.html} \quad ,

discovered (and proved) the following recurrences for $F_k(n)=F_k(n)(\alpha)$ for $1 \leq k \leq 10$.
They get increasingly complicated, and we will only state the first two in the body of this paper. From now on $a$ will stand for $\alpha$.

$$
-a \left(n +1\right) F_{1}\! \left(n \right)-\left(n + 1\right) F_{1}\! \left(n +1\right)+F_{1}\! \left(n +2\right) = 0 \quad .
\eqno(4)
$$
\vfill\eject
$$
4 a \left(2 n +5\right) \left(n +2\right) \left(n +1\right) \left(a +1\right)^{2} F_{2}\! \left(n \right)
$$
$$
+2 \left(n +2\right) \left(a +1\right) \left(4 a \,n^{2}+12 a n -4 n^{2}
+7 a -14 n -10\right) F_{2}\! \left(n +1\right)
$$
$$
-2 \left(n +2\right) \left(4 a n +4 n^{2}+8 a +16 n +17\right) F_{2}\! \left(n +2\right)+\left(2 n +3\right) F_{2}\! \left(n +3\right) = 0 \quad .
\eqno(5)
$$

For linear recurrences for $F_k(n)$ for $3 \leq k \leq 10$ see the output file

{\tt https://sites.math.rutgers.edu/\~{}zeilberg/tokhniot/oMourad1.txt} \quad .

These recurrences enable very fast computation of quite a few terms of these sequences. It is all implemented in the already mentioned Maple package {\tt Mouard.txt}. 
The direct url of this package is:

{\tt https://sites.math.rutgers.edu/\~{}zeilberg/tokhniot/Mourad.txt} \quad .

{\bf A User's Manual to the Maple package Mourad.txt}

To use it first download it to your favorite directory (usually Downloads). Start a Maple worksheet, make sure that the directory is the right one, and then type

{\tt read `Mouard.txt` ;} \quad .

To get a list of the procedures, type {\tt ezra();} . The main procedures are:

$\bullet$ {\tt Wder(L,a)}, that inputs a list of positive integers, $L=[k_1, \dots, k_n]$ and implements Equation $(2)$. This is very slow, and should not be used for large $L$.

$\bullet$ {\tt Operk(k,a,n,N)}, that inputs a positive integer $k$, and outputs the recurrence operator (where $N$ is the shift operator in $n$) annihilating
the sequence $F_k(n)$. These get very complicated for human eyes, but the computer does not mind and it enables a fast computation of many terms.

$\bullet$  {\tt SeqkF(k,K,a)}: inputs a positive integer $k$ and a positive integer $K$ and a symbol $a$ (what we called $\alpha$ above) and outputs the
first $K$ terms of the sequence of polynomials $F_k(n)=F_k(n)(a)$. For example, to get the weight enumerator of the set of permutations of a standard deck of cards where
no card can wind up at a location that originally was occupied by the same number ($1$ through $13$, where Jack, Queen, and King stand for $11,12,13$ respectively) 
but it is OK to have the same suit, type:

{\tt SeqkF(4,13,a)[13];} \quad ,

and you would get, in a few seconds, the following polynomial of degree $26$
$$
626486325682388256883179081695232\,a^{26} +
$$
$$
3948815860811007759557670403206807552\,a^{25}+
$$
$$
4226160446928101410675933447042193424384\,a^{24}+
$$
$$
1829313185198525509532452983498671376039936\,a^{23}+
$$
$$
425955227133577312273392421310068029118218240\,a^{22}+
$$
$$
61568711382255715699343414832865761752795578368\,a^{21}+
$$
$$
6015599331237497842549834616372527226200006852608\,a^{20}+
$$
$$
420030513102996289545618495318355347968579239673856\,a^{19}+
$$
$$
21779385529606788308065066752435641655566027030790144\,a^{18}+
$$
$$
861931009463580565142515454351924475556603802576486400\,a^{17}+
$$
$$
26556926811772603306934511893782498309330811792400580608\,a^{16}+
$$
$$
646219419386602045907824228576682527851206056554484727808\,a^{15}+
$$
$$
12544166147808400841334081628081554018739662394272604225536\,a^{14}+
$$
$$
195525408546538912690378251287680488219792943092212919435264\,a^{13}+
$$
$$
2455695605166443718371007842011087818790955115435503879454720\,a^{12}+
$$
$$
24867048146672227309345666989913796704728810126752820020379648\,a^{11}+
$$
$$
202569793911613274182929019082185092261014157201085369153486848\,a^{10}+
$$
$$
1320388339665569428585764027609539765653334771119656423470923776\,a^9+
$$
$$
6825167923093955037138102373992833000704975000443998456569135104\,a^8+
$$
$$
27602809328921835313793682068121303712142304270099611821308641280\,a^7+
$$
$$
85647342705993322148148235777401007447932223607159691210985046016\,a^6+
$$
$$
198159663830900044042641789039253865122617020230065397080602443776\,a^5+
$$
$$
327547473685724687587188995032714624999930689030717701980120154112\,a^4+
$$
$$
361148215004517312493645900517444844168859774724070502768740139008\,a^3+
$$
$$
234426065400514976953417524798811902707109969381695319447196139520\,a^2+
$$
$$
66394948050946830932484058263644488672722608355067055619597926400\,a \quad .
$$

Plugging-in $a=1$ and dividing by $4!^{13}$ (to get derangements where members of the same set are identified), we get
$$
1493804444499093354916284290188948031229880469556 \quad,
$$

agreeing with the title of [EkKZ], that handled the {\it straight}, rather than the {\it weighted}, enumeration of multiset derangements.

{\bf References}

[AlZ] Gert Almkvist and Doron Zeilberger, {\it The method of differentiating under the
integral sign}, J. Symbolic Computation {\bf 10}, 571-591 (1990). \hfill\break
{\tt https://sites.math.rutgers.edu/\~{}zeilberg/mamarim/mamarimhtml/duis.html} \quad .

[As] Richard Askey, {\it Orthogonal Polynomials and Special Functions}, SIAM, 1975.

[AsIs] R. Askey and M.E.H. Ismail, {\it Permutation problems and special functions}, Canad. J. Math. {\bf 28} (1976), 853-874.

[AsIsKo] R. Askey, M.E.H. Ismail, and T. Koornwinder, {\it Weighted permutation problems and Laguerre polynomials},
J. Comb. Theory (Ser. A.) {\bf 25} (1978), 277-287.

[AsIsRa] R. Askey, M.E.H. Ismail, and M.T. Rashed, {\it A derangement problem}. 
Technical report 1522, UW-Madison Mathematics Research Center, June 1975. 

[D] Robert Dougherty-Bliss, {\it Integral Recurrences from A to Z}, The American Mathematical Monthly {\bf 129} (2021), 805-815.

[EkKZ] Shalosh B. Ekhad, Christoph Koutschan, and Doron Zeilberger,
{\it There are EXACTLY $1493804444499093354916284290188948031229880469556$ ways to derange a standard deck of cards (ignoring suits) [and many other such useful facts]},
Enumerative Combinatorics and Applications {\bf 1} (2021), issue 3, article S2A17. \hfill\break
{\tt https://sites.math.rutgers.edu/\~{}zeilberg/mamarim/mamarimhtml/multider.html} \quad .

[EvG] S. Even and J. Gillis, {\it Derangements and Laguerre polynomials}, Math. Proc. Camb. Phil. Soc. {\bf 79} (1976), 135-143.

[FZ] Dominique Foata and Doron Zeilberger, {\it Laguerre polynomials, weighted derangements, and positivity},
SIAM J. Discrete Mathematics {\bf 1} (1988), 425-433.

[G] J. Gillis, {\it The statistics of derangement-a survey}, Journal of Statistical Physics {\bf 58} (1990),575-578.

[GIO]  J. Gillis, M. E. H. Ismail, and T. Offer, {\it An Asymptotic Problem in Derangement Theory},
SIAM Journal on Mathematical Analysis {\bf 21} (1990), 262-269.

[Z] Doron Zeilberger, {\it How Joe Gillis discovered combinatorial special function theory},
Math. Intell. {\bf 17(2)} (Spring 1995), 65-66. \hfill\break
{\tt https://sites.math.rutgers.edu/\~{}zeilberg/mamarim/mamarimhtml/gillis.html}

\bigskip
\hrule
\bigskip
Shalosh B. Ekhad, Department of Mathematics, Rutgers University (New Brunswick), Hill Center-Busch Campus, 110 Frelinghuysen
Rd., Piscataway, NJ 08854-8019, USA. \hfill\break
Email: {\tt ShaloshBEkhad at gmail  dot com}   \quad .
\bigskip
Doron Zeilberger, Department of Mathematics, Rutgers University (New Brunswick), Hill Center-Busch Campus, 110 Frelinghuysen
Rd., Piscataway, NJ 08854-8019, USA. \hfill\break
Email: {\tt DoronZeil at gmail  dot com}   \quad .

\bigskip

{\bf Jan. 7, 2025} 

\end